\documentclass{amsart}

\usepackage{amssymb}
\usepackage{times}
\usepackage{graphics}
\usepackage{array}
\usepackage{bbold}
\usepackage{stmaryrd}
\usepackage[all]{xy}

\newtheorem{thm}{Theorem}

\newtheorem{lem}[thm]{Lemma}

\theoremstyle{definition}

\newtheorem*{example}{Example}
\newtheorem*{ack}{Acknowledgement}

\theoremstyle{remark}

\newtheorem*{remark}{Remark}

\DeclareMathOperator{\tr}{trace}

\DeclareMathOperator{\End}{End}
\DeclareMathOperator{\Hom}{Hom}

\DeclareMathOperator{\Id}{Id}

\renewcommand{\k}{\mathbb{k}}

\newcommand{\onto}{\twoheadrightarrow}

\newcommand{\tto}{\longrightarrow}

\newcommand{\q}{\mathbf{q}}

\newcommand{\R}{\mathcal{R}}

\newcommand{\ZZ}{\mathbb{Z}}

\newcommand{\T}{\mathsf{T}}

\newcommand{\K}{\mathsf{K}}

\newcommand{\eA}{\underline{\rm end}\,A}
\newcommand{\sg}[1]{\underline{\rm end}\,#1}

\newcommand{\ul}[1]{\mathbf{#1}}

\begin{document}

%%%%%%%%%%%%%%%%%%% version %%%%%%%%%%%%%%%%%%%%%%%%%%%%%%%%%%%%%%%%%%%%
%%
%% \mbox{}
%%
%% \vspace{-1in}
%%
%% \hbox{}\hfill\parbox{2in}{\textsl{\small Version \#2 \\ \today \\
%%                          comments welcome}}
%% \vspace{.5in}
%%
%%%%%%%%%%%%%%%%%%%%%%%%%%%%%%%%%%%%%%%%%%%%%%%%%%%%%%%%%%%%%%%%%%%%%%%%

\title[MacMahon's Master Theorem]{Koszul algebras and the quantum MacMahon Master Theorem}

\author{Ph{\`u}ng H{\^o} Hai}
\address{Mathematik, University of Duisburg-Essen, Germany and
Institute of Mathematics, Hanoi, Vietnam}
\email{hai.phung@uni-duisburg-essen.de}
\thanks{PHH is supported by the DFG with a Heisenberg-Fellowship}

\author{Martin Lorenz}
\address{Department of Mathematics, Temple University,
    Philadelphia, PA 19122-6094}
\email{lorenz@temple.edu}
\thanks{ML's research is supported in part by NSA Grant H98230-05-1-0025 and
by Leverhulme Research Interchange Grant F/00158/X}


%\date{}


\begin{abstract}
We give a new proof of the quantum version of MacMahon's Master
Theorem due to Garoufalidis, L{\^e} and Zeilberger (one-parameter
case) and to Konvalinka and Pak (multiparameter case) by deriving it
from known facts about Koszul algebras.
\end{abstract}

\maketitle

%%%%%%%%%%%%%%%%%%%%%%%%%%%%%%%%%%%%%%%%%%%%%%%%%%%%%%%%%%%%%%%%%%%%%%%%

\section{Introduction}

%\subsection{} %\label{SS:}
In \cite{GLZxx}, Garoufalidis, L{\^e} and Zeilberger prove a quantum
version of MacMahon's celebrated ``Master Theorem"
\cite[pp.~97--98]{paMM60}.
%Percy Alexander MacMahon (26 Sept 1854 -- 25 Dec 1929)
As stated in \cite{GLZxx}, the generalization was motivated in part
by considerations in quantum topology and knot theory, and it also
answers a long-standing open question by G. Andrews \cite[Problem
5]{gA75}. An abundance of different proofs of the original Master
Theorem can be found in the literature; the quantum-generalization
in \cite{GLZxx} is proved by an application of the calculus of
difference operators developed by Zeilberger in \cite{dZ80}. Our
goal here is to derive the quantum MacMahon Master Theorem of
Garoufalidis, L{\^e}, and Zeilberger, along with its multiparameter
extension proved subsequently by Konvalinka and Pak \cite{mKiP},
fairly effortlessly from basic properties of Koszul algebras.
Indeed, the Koszul complex immediately leads to a generalized
MacMahon identity stated as equation \eqref{E:application} below.
The quantum MacMahon Master Theorem is the special case of
\eqref{E:application} where the Koszul algebra in question is the
so-called quantum affine space. Thus, neither the main result of
this note nor the methods employed are ours but we believe that the
connection between Koszul algebras and the quantum Master Theorem
deserves to be explicitly stated and fully exploited. The link to
quantum affine space and quantum matrices was in fact already
briefly mentioned in the last section of \cite{GLZxx}, but the proof
given here appears to be new.

We have tried to keep this note reasonably self-contained and
accessible to readers unfamiliar with Koszul algebras.
Sections~\ref{S:koszul} and \ref{S:bialgebras} serve to deploy the
pertinent background material concerning Koszul algebras and
characters in some detail. The operative technicalities for our
proof are collected in Lemmas~\ref{L:KA} and \ref{L:characters}
below; they are presented here with full proofs for lack of a
suitable reference. The quantum MacMahon Master Theorem
\cite[Theorem 1]{GLZxx}, \cite[Theorem 1.2]{mKiP} is then stated and
proved in Section~\ref{S:proof}. The short final
Section~\ref{S:modifying} discusses certain modifications of the
MacMahon identity.

Our basic reference for Koszul algebras are Manin's notes
\cite{yM88}; for bialgebras, our terminology follows Kassel
\cite{cK95}. We work over a commutative base field $\k$ except in
\ref{SS:grothendieck} and \ref{SS:characters} where $\k$ can be any
commutative ring at no extra cost. Throughout, $\otimes$ will stand
for $\otimes_{\k}$.

%\subsection{} %\label{SS:}
%blabla

%\begin{notat}
%Throughout, ...
%\end{notat}

%%%%%%%%%%%%%%%%%%%%%%%%%%%%%%%%%%%%%%%%%%%%%%%%%%%%%%%%%%%%%%%%%%%%%%%%

\section{Koszul algebras} \label{S:koszul}

\subsection{Quadratic algebras} \label{SS:quadratic}

A \emph{quadratic algebra} is a factor of the tensor algebra $\T(V)$
of some finite-dimen\-sional $\k$-vector space $V$ modulo the ideal
generated by some subspace $R(A) \subseteq \T(V)_2 = V^{\otimes 2}$.
Thus, $A \cong \T(V)/\left( R(A)\right)$.
%, where $\langle R(A)\rangle$ is the ideal of $\T(V)$ generated by $R(A)$.
The natural grading of $\T(V)$ descends to a grading $A =
\bigoplus_{d \ge 0} A_d$ of $A$ with $A_0 = \k$ and $A_1 \cong V$.
In practise, one often fixes a $\k$-basis
$\tilde{x}_1,\dots,\tilde{x}_n$ of $V$. Then $\T(V)$ can be viewed
as the free algebra $\k\langle\tilde{x}_1,\dots,\tilde{x}_n\rangle$.
The images $x_i = \tilde{x}_i \mod R(A)$ of the elements
$\tilde{x}_i$ in $A$ form a set of algebra generators for $A$.

\begin{example}[(Quantum affine $n$-space)]
This is the quadratic algebra $A = A_\q^{n|0}$ that is defined, for
a fixed scalars $0 \neq q_{ij} \in \k$ $(1 \le i < j \le n)$, as the
factor of $\k\langle\tilde{x}_1,\dots,\tilde{x}_n\rangle$ modulo the
ideal generated by the $2$-homogeneous elements
$\tilde{x}_j\tilde{x}_i - q_{ij}\tilde{x}_i\tilde{x}_j$ for $i < j$.
Thus, the algebra $A_\q^{n|0}$ is generated by elements
$x_1,\dots,x_n$ satisfying the relations
\begin{equation} \label{E:affinerel}
x_j x_i = q_{ij} x_i x_j\qquad (i < j).
\end{equation}
\end{example}

\subsection{Quadratic dual} \label{SS:dual}

Given a quadratic algebra $A \cong \T(V)/\left( R(A)\right)$,
consider the subspace $R(A)^\perp$ of the linear dual
$\left(V^{\otimes 2}\right)^*$ consisting of all linear forms on
$V^{\otimes 2}$ that vanish on $R(A)$. Identifying $\left(V^{\otimes
2}\right)^*$ with $V^{*\,\otimes 2}$ in the usual way, we may view
$R(A)^\perp \subseteq V^{*\,\otimes 2}$ and hence define a new
quadratic algebra by
$$
A^! = \T(V^*)/\left( R(A)^\perp\right) \ .
$$
The algebra $A^!$ is called the \emph{quadratic dual} of $A$. If
$\tilde{x}_1,\dots,\tilde{x}_n$ is a fixed $\k$-basis of $V$, as
above, then one usually chooses the dual basis
$\tilde{x}^1,\dots,\tilde{x}^n$ for $V^*$: $\langle \tilde{x}^i,
\tilde{x}_j\rangle = \delta_{i,j}$, where $\langle .\,, .\,\rangle$
denotes evaluation and $\delta_{i,j}$ is the Kronecker delta. This
yields algebra generators $x^i = \tilde{x}^i \mod R(A)^\perp$ for
$A^!$.

\begin{example}[(Quantum exterior algebra)]
Consider the algebra $A = A_\q^{n|0}$ as above. Following
\cite{yM88}, the quadratic dual $A^!$ will be denoted by
$A_\q^{0|n}$. The above procedure yields algebra generators
$x^1,\dots,x^n$ for $A_\q^{0|n}$ satisfying the defining relations
\begin{equation} \label{E:extrel1}
x^\ell x^\ell = 0
\end{equation}
for all $\ell$ and
\begin{equation} \label{E:extrel2}
x^k x^\ell + q_{k\ell} x^\ell x^k = 0 \qquad (k < \ell).
\end{equation}
%An element $r =
%\sum_{k,\ell} c_{k,\ell} \tilde{x}^k \otimes \tilde{x}^\ell \in
%(V^*)^{\otimes 2}$ belongs to $R(A)^\perp\rangle$ if and only if $0
%= \langle r,\tilde{x}_j\otimes\tilde{x}_i -
%q\tilde{x}_i\otimes\tilde{x}_j\rangle = c_{j,i} - q c_{i,j}$ holds
%for all $i < j$. Therefore, the elements $\tilde{x}^\ell \otimes
%\tilde{x}^\ell\ (\ell = 1,\dots,n)$ and $q\tilde{x}^j \otimes
%\tilde{x}^i + \tilde{x}^i \otimes \tilde{x}^j$ form a $\k$-basis of
%$R(A)^\perp\rangle$ and the quadratic dual
%$$
%A_q^{0|n} = \k\langle \tilde{x}^1,\dots,\tilde{x}^n \rangle/ \left(
%\tilde{x}^\ell\tilde{x}^\ell\,, \tilde{x}^i\tilde{x}^j + q
%\tilde{x}^j\tilde{x}^i \mid \ell, i < j \right) \ .
%$$
\end{example}

\subsection{The bialgebra $\eA$} \label{SS:eA}

Given quadratic algebras $A \cong \T(A_1)/\left( R(A)\right)$ and $B
\cong \T(B_1)/\left( R(B)\right)$, one defines the quadratic algebra
$A \bullet B = \T(A_1 \otimes B_1)/ \left( S_{23}(R(A) \otimes R(B))
\right)$. Here, $S_{23} \colon A_1^{\otimes 2} \otimes B_1^{\otimes
2} \to (A_1 \otimes B_1)^{\otimes 2}$ switches the second and third
factors.

The algebra $\eA = A^!\bullet A$ is particularly important.
Identifying $A_1$ with $V$ as above, we have
$$
\eA = \T(V^* \otimes V)/ \left( S_{23}(R(A)^\perp \otimes R(A))
\right) \ .
$$
Fix generators $x_i = \tilde{x}_i \mod R(A)$ for $A$ and $x^i =
\tilde{x}^i \mod R(A)^\perp$ for $A^!$ as in \S\ref{SS:quadratic}
and \S\ref{SS:dual}. Then the elements $\tilde{z}^j_i = \tilde{x}^j
\otimes \tilde{x}_i$ form a basis of $V^* \otimes V$ and their
images ${z}^j_i = \tilde{z}^j_i \mod R(\eA)$ form algebra generators
for $\eA$. The algebra $\eA$ is endowed with a comultiplication
\begin{equation} \label{E:Delta}
\Delta \colon \eA \to \eA \otimes \eA\ , \qquad \Delta(z_i^j) =
\sum_\ell z_i^\ell \otimes z_\ell^j
\end{equation}
and a counit
$$
\epsilon \colon \eA \to \k\ , \qquad \epsilon(z_i^j) = \delta_{i,j}
%\text{ (Kronecker delta)}
$$
which make $\eA$ a bialgebra over $\k$; see \cite[5.7 and
5.8]{yM88}. Furthermore, defining a coaction
\begin{equation} \label{E:dA}
\delta_A \colon A \to \eA \otimes A\ ,\qquad \delta_A(x_i) = \sum_j
z_i^j \otimes x_j\ ,
\end{equation}
the algebra $A$ becomes a left $\eA$-comodule algebra: $\delta_A$ is
a $\k$-algebra map that makes $A$ a comodule for $\eA$; see
\cite[5.4]{yM88}.

\begin{example}[(Right-quantum matrices)]
Returning to quantum affine space $A = A_\q^{n|0}$, we now describe
the defining relations between the generators ${z}^j_i$ $(i,j =
1,\dots,n)$ of the algebra $\sg{A_\q^{n|0}}$. Using the notation and
the relations of the Examples in \S\ref{SS:quadratic} and
\S\ref{SS:dual}, the foregoing leads to the following set of
generators for $R(\sg{A_\q^{n|0}}) \subseteq \left( V^* \otimes V
\right)^{\otimes 2}$: for all $\ell$ and $i < j$ we have a generator
$\tilde{z}^\ell_j \otimes \tilde{z}^\ell_i - q_{ij}\,
\tilde{z}^\ell_i \otimes \tilde{z}^\ell_j$\,, and for all $i < j$
and $k < \ell$ there is $\tilde{z}^k_j \otimes \tilde{z}^\ell_i -
q_{ij}\, \tilde{z}^k_i \otimes \tilde{z}^\ell_j + q_{k\ell}\,
\tilde{z}^\ell_j \otimes \tilde{z}^k_i - q_{k\ell}q_{ij}\,
\tilde{z}^\ell_i \otimes
\tilde{z}^k_j$. %After multiplying the second set of generators with
%$q^{-1}$,
Thus, we obtain the following relations between the generators
${z}^j_i$ of $\sg{A_\q^{n|0}}$\,:
\begin{equation}
z^\ell_j z^\ell_i = q_{ij}\, z^\ell_i z^\ell_j \quad\text{for all
$\ell$ and $i < j$}  \label{E:colrel}
\end{equation}
and
\begin{equation}
q_{ij}\, z^k_i z^\ell_j -  q_{k\ell}\, z^\ell_j z^k_i = z^k_j
z^\ell_i - q_{k\ell}q_{ij}\,\, z^\ell_i z^k_j \quad\text{for $i < j$
and $k < \ell$} \ .  \label{E:crossrel}
\end{equation}
The relations \eqref{E:colrel} and \eqref{E:crossrel} are called
\emph{column relations} and \emph{cross relations}, respectively;
they are identical with the relations considered in \cite[Section
1.2]{GLZxx} (for the one-parameter case) and in \cite[Section
1.3]{mKiP} Therefore, following the terminology of \cite{GLZxx},
\cite{mKiP}, we will call the $n \times n$-matrix $Z = (z_i^j)$ a
generic \emph{right-quantum $\q$-matrix}. Right-quantum
$\q$-matrices over an arbitrary $\k$-algebra $R$ are exactly the
matrices of the form $\varphi Z$ for an algebra map $\varphi \colon
\sg{A_\q^{n|0}} \to R$ (an ``$R$-point'' of the space defined by
$\sg{A_\q^{n|0}}$).
\end{example}

\subsection{The Koszul complex} \label{SS:KA}

For any quadratic algebra $A$, one can define Koszul complexes
$$
\K^{\ell,\bullet}(A) \colon \quad 0 \to A^{!\ *}_\ell \to A^{!\
*}_{\ell-1} \otimes A_1 \to \dots \to  A^{!\ *}_{1} \otimes
A_{\ell-1} \to A_\ell \to 0
$$
as in \cite[9.6]{yM88}. The quadratic algebra $A$ is called a
\emph{Koszul algebra} if all complexes $\K^{\ell,\bullet}(A)$ for
$\ell > 0$ are exact. It is known that if $A$ and $B$ are Koszul
then so are $A^!$ and $A \bullet B$. Moreover, if $A$ has a
so-called PBW-basis consisting of certain standard monomials, then
$A$ is Koszul; see, e.g., \cite[Theorem 4.3.1]{PP05}. This applies
in particular to quantum affine $n$-space $A = A_q^{n|0}$ which is
therefore Koszul \cite[4.2 Example 1]{PP05}.

\begin{lem} \label{L:KA}
Let $A$ be a quadratic algebra. Then all $A_i$ and all $
A^!_{j}{}^*$ are (left) comodules over $\eA$, and hence so are the
components $\K^{\ell,i}(A) =  A^{!\ *}_{\ell-i} \otimes A_i$ of the
Koszul complex. Moreover, the Koszul differential is an
$\eA$-comodule map.
\end{lem}

\begin{proof}
We will write $B = \eA$ for brevity; so $B_1 = V^* \otimes V$.
Equation \eqref{E:dA} shows that $\delta_A$ sends $V = A_1$ to $B_1
\otimes V$. Thus, $A_i$ is mapped to $B_i \otimes A_i \subseteq B
\otimes A_i$, and so each $A_i$ is a $B$-comodule. Moreover, as is
shown in \cite[5.5]{yM88}, the map $\delta_A$ comes from a map
$\tilde{\delta} \colon \T(V) \to \T(B_1) \otimes \T(V)$,
$\tilde{\delta}(\tilde{x}_i) = \sum_j \tilde{z}_i^j \otimes
\tilde{x}_j$, which satisfies $\tilde{\delta}(R(A)) \subseteq R(B)
\otimes \T(V) + \T(B_1) \otimes R(A)$. Following $\tilde{\delta}$ by
the canonical map $\T(B_1) \onto B = \T(B_1)/R(B)$ tensored with
$\Id_{\T(V)}$, we obtain a map
\begin{equation} \label{E:delta'}
\delta' \colon \T(V) \to B \otimes \T(V)\ , \quad
\delta'(\tilde{x}_i) = \sum_j z_i^j \otimes \tilde{x}_j
\end{equation}
satisfying $\delta'(V^{\otimes i}) \subseteq B_i \otimes V^{\otimes
i}$ and $\delta'(R(A)) \subseteq B \otimes R(A)$. Therefore, all
$V^{\otimes i}$ are $B$-comodules and $R(A)$ a $B$-subcomodule of
$V^{\otimes 2}$.
%It follows from the definition of
%$\eA$ that $R(A)$ is a left $\eA$-comodule, cf. \cite[5.5]{yM88}.
More generally, the subspaces $V^{\otimes i}\otimes R(A)\otimes
V^{\otimes m-2-i}$ for $0 \le i \le m-2$ are $B$-subcomodules of
$V^{\otimes m}$, and hence so are the subspaces
$\bigcap^{m-2}_{i=0}V^{\otimes i}\otimes R(A)\otimes V^{\otimes
m-2-i}$ and $R_m(A):= \sum^{m-2}_{i=0}V^{\otimes i}\otimes
R(A)\otimes V^{\otimes m-2-i}$. Following \cite[9.6]{yM88}, we
identify the former with the linear dual of $A^!_m = V^{*\,\otimes
m} / \sum^{m-2}_{i=0}V^{*\,\otimes i}\otimes R(A)^\perp \otimes
V^{*\,\otimes m-2-i}$; so
\begin{equation} \label{E:A!*}
A^{!\ *}_{m} = \bigcap^{m-2}_{i=0}V^{\otimes i}\otimes R(A)\otimes
V^{\otimes m-2-i} \ .
\end{equation}
Also according to  [\emph{ibid}, 9.6], the diagram
$$
\begin{array}{ccc}A^{!\ *}_{m+1} \otimes V^{\otimes n-1}&\hookrightarrow &A^{!\ *}_m \otimes V^{\otimes n}\\
\bigcup &&\bigcup\\
A^{!\ *}_{m+1} \otimes R_{n-1}(A)&\hookrightarrow & A^{!\ *}_m
\otimes R_n(A)
\end{array}
$$
induces the Koszul differential
$$
A^{!\ *}_{m+1} \otimes A_{n-1}=A^{!\ *}_{m+1} \otimes
\frac{V^{\otimes n-1}}{R_{n-1}(A)}\tto A^{!\ *}_{m} \otimes
\frac{V^{\otimes n}}{ R_n(A)}= A^{!\ *}_{m} \otimes A_{n}
$$
Since all spaces here are $B$-comodules, the differential is a
$B$-comodule map.
\end{proof}

%%%%%%%%%%%%%%%%%%%%%%%%%%%%%%%%%%%%%%%%%%%%%%%%%%%%%%%%%%%%%%%%%%%%%%%%

\section{Bialgebras and characters} \label{S:bialgebras}

\subsection{The Grothendieck ring} \label{SS:grothendieck}

For now, let $B$ denote an arbitrary bialgebra over some commutative
base ring $\k$. We let $\R_B$ denote the Grothendieck ring of all
(left) $B$-comodules that are finitely generated (f.g.) projective
over $\k$. Thus, for each such $B$-comodule, $V$, there is an
element $[V] \in \R_B$ and any short exact sequence $0 \to U \to V
\to W \to 0$ of $B$-comodules (f.g.~projective over $\k$) gives rise
to an equation $[V] = [U] + [W]$ in $\R_B$. Multiplication in $\R_B$
is given by the tensor product of $B$-comodules; see
%\cite[1.8.2]{sM93}
\cite[III.6]{cK95}.

\subsection{Characters} \label{SS:characters}

Continuing with the notation of \S\ref{SS:grothendieck}, let $V$ be
a $B$-comodule that is f.g.~projective over $\k$. The structure map
$\delta_V \colon V \to B \otimes V$ is an element of
$\Hom_\k(V,B\otimes V)$. Using the standard isomorphisms
$\Hom_\k(V,B\otimes V) \cong B \otimes V \otimes V^*$ (see, e.g.,
\cite[II.4.2]{nB70}) and letting $\langle .\,, .\,\rangle \colon V
\otimes V^* \to \k$ denote the evaluation map, we have a
homomorphism
$$
\xymatrix{ \Hom_\k(V,B\otimes V) \cong B \otimes V \otimes V^*\quad
%\stackrel{\langle .\,, .\,\rangle \otimes \Id_B}{\tto}
\ar[r]^-{\Id_B \otimes \langle .\,, .\,\rangle} & \quad  B \otimes
\k \cong B \ . }
$$
The image of $\delta_V$ under this map will be denoted by $\chi_V$.
If $V$ is free over $\k$, with basis $\{v_i\}$, and $\delta_V(v_j) =
\sum_i b_{i,j} \otimes v_i$ then
\begin{equation}
\chi_V = \sum_i b_{i,i}\ . \label{E:trace}
\end{equation}

\begin{lem} \label{L:characters}
The map $[V] \mapsto \chi_V$ yields a well-defined ring homomorphism
$\chi \colon \R_B \to B$.
\end{lem}

\begin{proof}
The assertion is equivalent to the identities $\chi_{V \otimes W} =
\chi_{V}\chi_{W}$ for any two $B$-comodule $V$ and $W$
(f.g.~projective over $\k$) and $\chi_V = \chi_U + \chi_W$ for any
short exact sequence $0 \to U \stackrel{\mu}{\tto} V
\stackrel{\pi}{\tto} W \to 0$. Both are easy to check; we sketch the
proof of the second identity. Fix $\k$-linear splittings $\sigma
\colon W \to V$ and $\tau \colon V \to U$ with $\tau \circ \mu =
\Id_U$, $\pi \circ \sigma = \Id_W$ and $\sigma \circ \pi + \mu \circ
\tau = \Id_V$. Under the canonical isomorphism $V \otimes V^* \cong
\End_{\k}(V)$, evaluation $\langle .\,, .\,\rangle$ becomes the
trace map $\tr_V \colon \End_{\k}(V) \to \k$. For any $\varphi \in
\End_{\k}(V)$, we have $\tr_V(\varphi) = \tr_U(\tau \circ \varphi
\circ \mu) + \tr_W(\pi \circ \varphi \circ \sigma)$.
%Writing $\delta_V = \sum_i b_i \otimes \varphi_i  \in B \otimes \End_{\k}(V)$,
The equation $\delta_V \circ \mu = (\Id_B \otimes \mu) \circ
\delta_U$ gives $\delta_U = (\Id_B \otimes \tau) \circ \delta_V
\circ \mu$. Similarly, $\delta_W = (\Id_B \otimes \pi) \circ
\delta_V \circ \sigma$. Therefore, $\chi_V = (\Id_B \otimes
\tr_V)(\delta_V) = (\Id_B \otimes \tr_U)(\delta_U) + (\Id_B \otimes
\tr_W)(\delta_W) =  \chi_U + \chi_W$, as desired.
\end{proof}

\subsection{Application to Koszul algebras} \label{SS:application}

Returning to the case of a Koszul algebra $A$ over a field $\k$, we
apply the foregoing to the bialgebra $B = \eA$. By Lemma~\ref{L:KA},
the (exact) Koszul complex $\K^{\ell,\bullet}(A)$ for $\ell > 0$
gives an equation in $\R_{B}$\,:
$$
\sum_i (-1)^i [ A^{!\ *}_{i}][A_{\ell-i}] = 0 \ .
$$
Equivalently, defining the Poincar{\'e} series $P_A(t) = \sum_i
[A_{i}] t^i$ and $P_{A^{\!!*}}(t) = \sum_i [A^{!\ *}_{i}] t^i$, we
have
$$
P_A(t) P_{A^{\!!*}}(-t) = 1
$$
in the power series ring $\R_{B}\llbracket t \rrbracket$. Applying
the ring homomorphism $\chi\llbracket t \rrbracket \colon
\R_{B}\llbracket t \rrbracket \to B\llbracket t \rrbracket$
(Lemma~\ref{L:characters}), this equation becomes the following
equation in $B\llbracket t \rrbracket$:
\begin{equation} \label{E:application}
\left( \sum_{\ell \ge 0} \chi_{A_\ell} t^\ell \right) \cdot \left(
\sum_{m \ge 0} (-1)^m \chi_{A^{!\,*}_m} t^m \right) = 1 \ .
\end{equation}
Since the coactions $\delta_{A_\ell}$ and $\delta_{A^{!\,*}_m}$ in
\eqref{E:dA} and \eqref{E:delta'} respect the grading, both factors
actually belong to the subring $\prod_{n\ge 0} B_n t^n$ of
$B\llbracket t \rrbracket$.

%%%%%%%%%%%%%%%%%%%%%%%%%%%%%%%%%%%%%%%%%%%%%%%%%%%%%%%%%%%%%%%%%%%%%%%%

\section{Proof of the quantum MacMahon Master Theorem} \label{S:proof}

The quantum MacMahon Master Theorem \cite[Theorem 1]{GLZxx},
\cite[Theorem 1.2]{mKiP} is the special case of
equation~\eqref{E:application} where $A$ is quantum affine space. We
need to evaluate the characters $\chi_{A_\ell}$ and
$\chi_{A^{!\,*}_m}$ for $A = A_\q^{n|0}$.

Choose generators $x_1,\dots,x_n$ for $A$ as in
\S\ref{SS:quadratic}. For each $n$-tuple $\ul{m} =
(m_1,\dots,m_n)\in \ZZ_{\geq 0}^n$, put
$x^{\ul{m}}:=x_1^{m_1}x_2^{m_2}\ldots x_n^{m_n} \in A$. Then the
homogeneous component $A_\ell$ has a $\k$-basis consisting of the
elements $x^{\ul{m}}$ with $|\ul{m}| = m_1+m_2+\ldots+m_n = \ell$.
(This is the PBW-basis of $A_\q^{n|0}$ that was referred to
earlier.) With respect to this basis, the coaction $\delta_A$ of
$B=\eA$ on $A_\ell$ in \eqref{E:dA} has the form
$$
\delta_A(x^{\ul{m}}) = \delta_A(x_1)^{m_1} \delta_A(x_2)^{m_2}
\cdots \delta_A(x_n)^{m_n} = \sum_{\substack{\ul{r} \in \ZZ_{\geq
0}^n
\\ |\ul{r}| = \ell}} b_{\ul{r},\ul{m}} \otimes x^{\ul{r}}
$$
for uniquely determined $b_{\ul{r},\ul{m}} \in B_\ell$. In
particular, $G(\ul{m}):= b_{\ul{m},\ul{m}} \in B_{|\ul{m}|}$ has the
same meaning as in \cite{GLZxx} and \cite{mKiP}. From equation
\eqref{E:trace} we obtain the formula
\begin{equation} \label{E:chiA}
\chi_{A_\ell} = \sum_{\ul{m} \colon |\ul{m}| = \ell} G(\ul{m}) \ .
\end{equation}

It remains to calculate the character $\chi_{A^{!\,*}_m}$. To this
end, we identify $A^{!\ *}_m$ with the subspace $R_m(A^{!})^\perp$
of $\T(V)_m = V^{\otimes m}$ as in \eqref{E:A!*} and we think of
$\T(V)$ as the free algebra
$\k\langle\tilde{x}_1,\dots,\tilde{x}_n\rangle$ as in
\S\ref{SS:quadratic}.  Then $A^{!\ *}_m$ has a $\k$-basis consisting
of the elements
$$
\wedge \tilde{x}_{J}:= \sum_{\pi\in\frak S_m} w(\pi)
\tilde{x}_{j_{\pi1}} \otimes \tilde{x}_{j_{\pi2}} \otimes \ldots
\otimes \tilde{x}_{j_{\pi m}}\ ,
$$
where $J=(j_1<j_2<\ldots< j_m)$ is an $m$-element subset of
$\{1,\dots,n\}$, $\frak S_m$ is the symmetric group on
$\{1,\dots,m\}$ and
$$
w(\pi) = \prod_{i < j,\, \pi i
> \pi j} \left(-q_{\pi j,\pi i}\right)^{-1}
$$
as in \cite{mKiP}. Indeed, using the generators of $R(A^!)$
exhibited in \eqref{E:extrel1} and  \eqref{E:extrel2}, it is
straightforward to check that $\wedge \tilde{x}_{J}$ vanishes on
$R_m(A^!) = R(A^!) \cap T(V^*)_m$. Hence, $\wedge \tilde{x}_{J}$
belongs to $A^{!\ *}_m$. Furthermore, the elements $\wedge
\tilde{x}_{J}$ for distinct $m$-element subsets $J \subseteq
\{1,\dots,n\}$ are obviously $\k$-linearly independent and their
number is $n\choose m$ which is equal to the dimension of $A^!_m$.
Therefore, we obtain the claimed basis of $A^{!\ *}_m$.

Consider the basis vector $\wedge \tilde{x}_{[1,n]}$ corresponding
to the subset $J=\{1,2,...,n\} = :[1,n]$.  The coaction $\delta'$ of
$B = \eA$ on $A^{!\ *}_m$ in equation \eqref{E:delta'} satisfies
\begin{equation} \label{E:qdet'}
\delta'(\wedge \tilde{x}_{[1,n]})= {\det}_\q(Z) \otimes
\wedge\tilde{x}_{[1,n]} \ ,
\end{equation}
where $Z = (z_i^j)$ is the generic $n \times n$ right-quantum
$\q$-matrix as in \S\ref{SS:eA} and
\begin{equation} \label{E:qdet}
{\det}_\q(Z) = \sum_{\pi\in\frak S_m}w(\pi) z^{1}_{\pi1}
z^{2}_{\pi2} \dots z^{n}_{\pi n}
\end{equation}
denotes the multiparameter quantum determinant of $Z$ as defined in
\cite{AST91} (see also \cite{mKiP}). To prove \eqref{E:qdet'}, note
that the element $\wedge\tilde x_{[1,n]}$ spans the one-dimensional
$B$-subcomodule $A^{!\ *}_n$ of $V^{\otimes n}$. Hence, we certainly
have
\begin{equation} \label{E:qdet''}
\delta'(\wedge \tilde{x}_{[1,n]})= D \otimes \wedge\tilde{x}_{[1,n]}
\end{equation}
for some group-like element $D \in B$; cf. \cite[8.2]{yM88}. We need
to show that $D = \det_\q(Z)$. But equation \eqref{E:delta'} readily
implies that
$$
\delta'(\tilde{x}_{\pi1}\otimes \ldots \otimes \tilde{x}_{\pi n}) =
z^1_{\pi 1}\ldots z^n_{\pi n} \otimes \tilde{x}_1
\otimes\ldots\otimes \tilde{x}_n + \text{ other terms} \ ,
$$
where each of the ``other terms" belongs to some $B \otimes
\tilde{x}_{i_1} \otimes\ldots\otimes \tilde{x}_{i_n}$ with
$(i_1,\dots,i_n) \neq (1,\dots,n)$. Multiplying with $w(\pi)$ and
summing up over $\pi$ we see that $\delta'(\wedge \tilde{x}_{[1,n]})
= \det_\q(Z) \otimes \tilde{x}_1 \otimes\ldots\otimes \tilde{x}_n +
\text{ other terms}$. In view of \eqref{E:qdet''}, this implies
\eqref{E:qdet'}.

For a general $m$-element subset $J = (j_1<j_2<\ldots< j_m)
\subseteq \{1,\dots,n\}$, let $Z_J$ denote the submatrix
$(z_i^j)_{i,j \in J}$ of $Z$; this is also a generic right-quantum
matrix. Equation \eqref{E:qdet'} implies that $ \delta'(\wedge
\tilde{x}_{J}) = {\det}_\q(Z_J)\otimes \wedge \tilde{x}_{J} + \rho$,
where $\rho \in \sum_{I\neq J} B \otimes \wedge \tilde{x}_{I}$ and
\begin{equation} \label{E:qdetZJ}
{\det}_\q(Z_J) = \sum_{\pi\in\frak S_m}w(\pi) z^{j_1}_{j_{\pi1}}
z^{j_2}_{j_{\pi2}} \dots z^{j_m}_{j_{\pi m}} \ .
\end{equation}
By \eqref{E:trace}, the character of $A^{\!!\,*}_m$ is therefore
given by
\begin{equation} \label{E:chiA!*}
\chi_{A^{\!!\,*}_m} = \sum_{\substack{J \subseteq \{1,\dots,n\}
\\ |J| = m}} {\det}_\q(Z_J) \ .
\end{equation}

To summarize, we rephrase equation~\eqref{E:application} with
\eqref{E:chiA} and \eqref{E:chiA!*} using the notation of
\cite{GLZxx}:

\begin{thm}[(\cite{GLZxx}, \cite{mKiP})] \label{T:GLZ}
Let $A$ denote the $\k$-algebra generated by $x_1, \dots, x_n$
subject to the relations \eqref{E:affinerel}, let $B$ be the
$\k$-algebra generated by $z_1^1, \dots, z_n^n$ subject to
\eqref{E:colrel} and \eqref{E:crossrel}, and let $Z = (z_i^j)_{n
\times n}$ denote the generic right-quantum $\q$-matrix. Consider
the elements $X_i = \sum_j z_i^j \otimes x_j \in B \otimes A =
\bigoplus_{\ul{m}} B \otimes x^{\ul{m}}$, where $\ul{m}$ runs over
the $n$-tuples $\ul{m} = (m_1,\dots,m_n)\in \ZZ_{\geq 0}^n$ and
$x^{\ul{m}}:=x_1^{m_1}x_2^{m_2}\ldots x_n^{m_n}$. Define power
series in $B\llbracket t \rrbracket$ by
$$
\mbox{\rm Bos}(Z):= \sum_{\ell \ge 0}\ \sum_{\ul{m} \colon |\ul{m}|
= \ell} G(\ul{m}) t^\ell  \ ,
$$
where $G(\ul{m})$ is the $B$-coefficient of $x^{\ul{m}}$ in
$X^{\ul{m}} = X_1^{m_1}X_2^{m_2}\ldots X_n^{m_n}$, and
$$
\mbox{\rm Ferm}(Z):= \sum_{m \ge 0}\ \sum_{\substack{J \subseteq
\{1,\dots,n\}
\\ |J| = m}} (-1)^m {\det}_\q(Z_J) t^m
$$
with ${\det}_\q(Z_J)$ as in \eqref{E:qdetZJ}. Then:
$$
\mbox{\rm Bos}(Z)\cdot  \mbox{\rm Ferm}(Z) = 1 \ .
$$
\end{thm}

%%%%%%%%%%%%%%%%%%%%%%%%%%%%%%%%%%%%%%%%%%%%%%%%%%%%%%%%%%%%%%%%%%%%%%%%

\section{Modifying the MacMahon identity} \label{S:modifying}

Applying endomorphisms of $B\llbracket t \rrbracket$ to the
generalized MacMahon identity \eqref{E:application}, we obtain new
versions of this identity. In this section, we discuss a particular
example for the case $A = A_\q^{n|0}$. As usual, we let $B = \eA$.

The algebraic torus $(T \times T)/\k^*$, with $T = (\k^*)^n$ and
$\k^*$ diagonally embedded in $T \times T$, acts on $B$ (and hence
on $B\llbracket t \rrbracket$) via
\begin{equation*} \label{E:tau}
(\tau,\tau')(z_i^j)= c_i d_j^{-1} z_i^j
\end{equation*}
for $\tau = (c_1,\dots,c_n), \tau' = (d_1,\dots,d_n) \in  T$.
Indeed, $(\tau,\tau')$ respects the relations \eqref{E:colrel} and
\eqref{E:crossrel} of $B$, and hence $(\tau,\tau')$ defines a graded
$\k$-algebra automorphism of $B$. Note that the diagonal subgroup
$\{(\tau,\tau) \mid \tau \in T \}$ acts by bialgebra automorphisms
while the subgroups $T \times \{1\}$ and $\{1\} \times T$ act by
left and right $B$-comodule automorphisms, respectively. For any
subset $J \subseteq \{1,\dots,n\}$, we have
\begin{equation*}
(\tau,\tau')({\det}_\q(Z_J)) = \mu_J(\tau \tau'^{-1}) {\det}_\q(Z_J)
\quad\text{with $\mu_J(c_1,\dots,c_n) =  \prod_{j \in J} c_j$}\ ,
\end{equation*}
and for $\ul{m} = (m_1,\dots,m_n)\in \ZZ_{\geq 0}^n$,
\begin{equation*}
(\tau,\tau')(G(\mathbf m)) = \mu_{\ul{m}}(\tau \tau'^{-1}) G(\mathbf
m) \quad\text{with $\mu_{\ul{m}}(c_1,\dots,c_n) = \prod_i
c_i^{m_i}$.}
\end{equation*}
Therefore, applying $(\tau, \tau')$ to the characters of $A_\ell$
and $A^{!\ *}_m$ as determined in \eqref{E:chiA} and
\eqref{E:chiA!*}, we obtain
\begin{eqnarray}
(\tau, \tau')(\chi_{A^{\!!\,*}_m})&=&\sum_{\substack{J \subseteq
\{1,\dots,n\} \\ |J| = m}}
\mu_J(\tau \tau'^{-1}){\det}_\q(Z_J)  \\
(\tau, \tau')(\chi_{A_\ell})&=&\sum_{|\mathbf m|=\ell}
\mu_{\ul{m}}(\tau \tau'^{-1}) G(\mathbf m)\ .
\end{eqnarray}
The particular choice
\begin{equation}\label{E:special}
\tau = (q^{n-1}, q^{n-3},\ldots, q^{1-n})
\end{equation}
and $\tau' =1$ leads to the following version of quantum MacMahon
Master theorem:
\begin{equation}\label{eq6}
\widetilde{\text{Bos}}(Z)\cdot  \widetilde{\text{Ferm}}(Z) = 1 \ ,
\end{equation}
where
\begin{eqnarray}
\widetilde{\text{Bos}}(Z)&:=& \sum_{\ell \ge 0}\ \sum_{\ul{m} \colon
|\ul{m}|
= \ell}q^{\ell(n+1)-2\sum im_i} G(\ul{m}) t^\ell  \\
\widetilde{\text{Ferm}}(Z)&:=& \sum_{m \ge 0}\ \sum_{J = (j_1 <
\dots < j_m)} (-1)^mq^{m(n+1)-2(j_1+j_2+\ldots+j_m)} {\det}_\q(Z_J)
t^m \ .
\end{eqnarray}

\begin{remark}
Let $H:=\underline{\text{gl}}(A)$ be the coordinate ring on the
quantum general linear group (cf. \cite[8.5]{yM88}). Then $H$ is a
Hopf algebra which is in fact coquasitriangular or cobraided (cf.
\cite[VIII.5]{cK95} ). Thus for any finite dimensional comodule $X$,
there exists a canonical isomorphism $X\to X^{**}$ of $H$-comodules,
given in terms of the braiding. This isomorphism is in general not
compatible with the tensor product. We note that the category of
$H$-comodules also possesses a ribbon \cite[XIV.6]{cK95}.
%[ibid,XIV.6].
By composing the above canonical isomorphism with the ribbon, one
obtains a functorial isomorphism $\tau_X:X\to X^{**}$ which is
compatible with the tensor product in the sense that
$$\tau_{X\otimes Y}=\tau_X\otimes \tau_Y \ .$$
Using $\tau_X$ we  can defined a new type of character of $X$,
called quantum character, as follows (cf. \cite{phH02}).
$\chi_{q,X}$ is the image of $1\in \k$ under the map
$$\k \stackrel{\text{db}}\tto X^*\otimes X\stackrel{\Id\otimes\tau_X}\tto X^*\otimes X^{**}\tto
 H\otimes X^*\otimes X^{**}\stackrel{\text{ev}}\tto
 H$$
 where $\text{db}$ is the ``dual base" map  $1\mapsto \sum_i x^i\otimes x_i$.
As for the ordinary character, one can show that the quantum
character is multiplicative with respect to the tensor product and
additive with respect to exact sequences. Applying the quantum
character to the Koszul complex in \S\ref{SS:KA} we obtain an
identity analogous to \eqref{E:application} for $\chi_{q,A_\ell}$
and $\chi_{q,A_m^{!*}}$.

 Explicit computation shows that, for $X=V=A_1$, $\tau_{V}=\text{diag}(q^{n-1},
 q^{n-3},\ldots, q^{1-n})$. Further, one can check that
 $$\chi_{q,A_\ell}=\tau(\chi_{A_\ell});\quad \chi_{q,A_m^{!*}}=
 \tau(\chi_{A_m^{!*}})$$
This explains the origin of the choice of $\tau$ in
\eqref{E:special}.
\end{remark}

%%%%%%%%%%%%%%%%%%%%%%%%%%%%%%%%%%%%%%%%%%%%%%%%%%%%%%%%%%%%%%%%%%%%%%%%

\begin{ack}
Part of the work on the final version of this article was done while
ML was visiting the Universit{\' e} de Montpellier in June 2006. It
is a pleasure to thank Claude Cibils for arranging the visit and for
his warm hospitality. Moreover, ML wishes to thank S.~Paul Smith for
providing him with a copy of \cite{yM88}. Finally, the authors are
grateful to Roland Berger for his suggestion to include the
multiparameter case of the Master Theorem; the original version of
this article dealt with the one-parameter case only.
\end{ack}

%%%%%%%%%%%%%%%%%%%%%%%%%%%%%%%%%%%%%%%%%%%%%%%%%%%%%%%%%%%%%%%%%%%%%%%%

\def\cprime{$'$}
\providecommand{\bysame}{\leavevmode\hbox
to3em{\hrulefill}\thinspace}
\providecommand{\MR}{\relax\ifhmode\unskip\space\fi MR }
\providecommand{\MRhref}[2]{%
  \href{http://www.ams.org/mathscinet-getitem?mr=#1}{#2}
} \providecommand{\href}[2]{#2}

%%%%%%%%%%%%%%%%%%%%%%%%%%%%%%%%%%%%%%%%%%%%%%%%%%%%%%%%%%%%%%%%%%%%%%%%

\end{document}